\documentclass[12pt]{amsart}
\usepackage{amsthm}
\usepackage{amsmath}
\usepackage{amsfonts} 

\newtheorem{theorem}{Theorem}[section]
\newtheorem{lemma}[theorem]{Lemma}
\newtheorem{proposition}[theorem]{Proposition}
\newtheorem{corollary}[theorem]{Corollary}
\theoremstyle{definition}

\theoremstyle{remark}
\newtheorem{remark}[theorem]{Remark}
\newcommand{\Rd}{{\rm Res}_{d_{1},\dots ,d_{n}}}

\newcommand{\efen}{f_1,\dots,f_n}
\newcommand{\im}{{\rm Im}}
\newcommand{\B}{{\rm Bez}}
\newcommand{\ba}{{\mathcal B}}
\newcommand{\ca}{{\mathcal C}}
\def\EE{{\mathbb E}}

\def\NN{{\mathbb N}}
\def\ZZ{{\mathbb Z}}

\def\CP{{\mathbb P}}

\def\F{{\mathcal F}}

\newenvironment{con}{\noindent{\it Convention}.}{}

\begin{document}

\title{Explicit Formulas for the Multivariate 
Resultant }
\author{Carlos D' Andrea and Alicia Dickenstein}
\thanks{Both authors are supported by Universidad de Buenos Aires, grant TX094}
\thanks{The second author is also supported by CONICET, Argentina, 
and the Wenner-Gren Foundation, Sweden.}
\date{}

\begin{abstract}
We present formulas for  the multivariate resultant
as a quotient of two determinants. They extend the classical Macaulay 
formulas, and involve matrices of considerably smaller size, whose non
zero entries 
include coefficients of the given polynomials
and coefficients of their Bezoutian.
These formulas can also be viewed as an explicit computation of
the morphisms and the determinant of a resultant complex.
\end{abstract}

\maketitle

\section{Introduction}

Given $n$ homogeneous polynomials $f_1,\dots,f_n$ in $n$ variables
over an algebraically closed field $k$
with respective degrees $d_1,\dots, d_n$, the resultant $\Rd(\efen)$ 
is an irreducible
polynomial in the coefficients of $\efen$, which
vanishes whenever $\efen$ have a common root in projective space.
The study of resultants goes back to classical work of Sylvester, B\'ezout,
Cayley, Macaulay and Dixon. The use of resultants as a computational tool 
for elimination of variables as well as a tool for the study of complexity 
aspects of polynomial system solving in the last decade,
 has renewed the interest in finding explicit formulas for their computation
(cf. \cite{AS}, \cite{C1}, \cite{C2}, \cite{cm2},\cite{EM}, \cite{kps},\cite{L}, 
\cite{R},\cite{rojas}). 
\par
By a determinantal formula it is meant a matrix 
whose  entries are polynomials in the coefficients of $\efen$
and whose determinant equals the
resultant $\Rd(\efen)$. Of course, the interest on such a formula
is the computation of the resultant, and so it is implicit that
the entries should be algorithmically computed from the
inputs. It is also meant that all non-zero entries have degree strictly
less than the degree of the resultant.
\par
In case all $d_i$ have a common value $d$, all currently
known determinantal formulas
are listed by Weyman and Zelevinsky in \cite{wz}.
This list is short:  if  $d \geq 2,$ 
there exist determinantal formulas for all $d$ just
for binary forms (given by the well known Sylvester matrix),  
ternary forms and quaternary
forms; when $n=5$, the only possible values for $d$ are $2$ and $3$; finally, for
$n=6$, there exists a determinantal formula only for $d=2$.
We find similar strict restrictions  on general $n, d_1,\dots,
d_n$ (cf. Lemma \ref{square}).
\par
Given $d_1,\dots,d_n,$ denote $t_n:= \sum_{i=1}^n (d_i-1)$ the critical degree. 
Classical Macaulay formulas \cite{Mac}
describe the resultant $\Rd(\efen)$
as an explicit quotient of two determinants. 
These formulas involve a matrix of size at least the number of monomials
in $n$ variables of degree $t_n+1$, and a submatrix of it. 
\par
Macaulay's work has been revisited and sharpened by Jouanolou 
in \cite{jou}, where he proposes for each 
$t\geq0,$  a square matrix $M_t$ of size
\begin{equation}
\label{tamatrix}
\rho\left(t\right):=\binom{t+n-1}{n-1}+ i(t_n-t)
\end{equation}
whose determinant is a nontrivial multiple of $\Rd(\efen)$ 
(cf. \cite{jou}, $3.11.19.7$).
Here, $i(t_n-t)$ denotes the dimension of the $k-$vector space
of elements of degree $t_n-t$ in the ideal generated by a
regular sequence  of $n$ polynomials with degrees $d_1,\dots,d_n.$
Moreover, Jouanolou shows that 
the resultant may be computed as the ratio between the determinant
 of $M_{t_n}$ and
the determinant of one of its square submatrices.
(cf. \cite{jou}, Corollaire $3.9.7.7$).  

\smallskip

In this paper, we explicitly find  the extraneous factor in Jouanolou's 
formulation, i.e. the polynomial $\det(M_t) /\Rd(\efen)$, for all 
$t\geq 0$  which again happens
to be the determinant of a submatrix ${\EE}_t$ of $M_t$ for every $t$, 
and this allows us to present new resultant formulas {\it \`a la 
Ma\-cau\-lay \/} for the resultant, i.e. as a quotient of two determinants
\begin{equation}
\label{macaulay}
 \Rd(\efen) = \frac{\det(M_t)} {\det({\EE}_t)}.
\end{equation}
For $t>t_n,$  we recover 
Macaulay's classical formulas.  For $t \leq t_n,$ the size of the matrix
$M_t$ is considerably smaller.  

\smallskip

In order to give explicit examples, we need to recall the definition
of the {\it Bezoutian associated with $\efen$}  
(cf. \cite{BCRS}, \cite{fgs}, \cite{kunz}, \cite{ss} and  \cite{jou} under
the name   ``Formes de Morley''). 
Let $\left(f_1,\dots,f_n\right)$ be a sequence of generic homogeneous polynomials
with respective degrees $d_1,\dots,d_n$ 
\begin{equation*}
f_{i}:=\sum_{|\alpha _{i}|=d_{i}}a_{\alpha _{i}}X^{\alpha _{i}} \quad
\in A\left[ X_{1},\dots ,X_{n}\right] ,
\end{equation*}
where 
$A$ is the factorial domain $A:={\ZZ}\left[ a_{\alpha _{i}}
\right]_{|\alpha _{i}|= d_{i},i=1,\dots ,n}.$
\par
Introduce two sets of $n$ variables $X,Y$  
and for each pair $(i,j)$ with  $ 1 \leq i,j \leq n$, 
write $\Delta_{ij}(X,Y)$ for the incremental quotient
\begin{equation}
\label{deltaij}
\frac{f_{i}(Y_{1},\dots ,Y_{j-1},X_{j},\dots
,X_{n})-f_{i}(Y_{1},\dots ,Y_{j},X_{j+1},\dots ,X_{n})}{X_{j}-Y_{j}}.
\end{equation}
Note that $f_{i}(X)-f_{i}(Y)=\sum_{j=1}^{n}\Delta _{ij}(X,Y)(X_{j}-Y_{j}).$

The determinant
\begin{equation}
\label{delta}
\Delta(X,Y):= \det(\Delta_{ij}(X,Y))_{1\leq i,j\leq n}=\sum_{\left| \gamma 
\right| \leq t_{n}}\Delta_{\gamma }\left( X\right) .Y^{\gamma }.
\end{equation}
is a representative of the
{\it Bezoutian} associated with 
$\left(f_1,\dots,f_n\right).$  
It is a homogeneous polynomial in $A \left[ X,Y\right] $  
of degree $t_{n}.$ 
\par
Recall also that
$$\deg \Rd(\efen) = \sum_{i=1}^n d_1 \dots d_{i-1} \cdot d_i \dots d_n.$$

\smallskip

As a first example, let  $n=3, \ \left(d_1,d_2,d_3\right)=\left(1,1,2\right),$
and let
$$
\begin{array}{ccc}
f_1 &= &a_1 X_1 + a_2 X_2 + a_3 X_3\\
f_2 &=& b_1 X_1 + b_2 X_2 + b_3 X_3\\
f_3 &=& c_1 X_1^2 + c_2 X_2^2 + c_3 X_3^2+ c_4 X_1X_2 + c_5 X_1 X_3 + c_6 X_2
X_3\\
\end{array}
$$
be generic polynomials of respective degrees $1,1,2$. 
Here, $t_3= 1.$
Macaulay's classical matrix $M_2$ looks as follows:
$$\left(
\begin{array}{cccccc}
a_1 & 0 & 0 & 0 & 0 & c_1 \\
0 & a_2 & 0 & b_2 & 0 & c_2 \\
0 & 0 & a_3 & 0 &b_3 & c_3 \\
a_2 & a_1 & 0 & b_1 & 0 & c_4 \\
a_3 & 0 & a_1 & 0 & b_1 & c_5 \\
0 &a_3& a_2 & b_3 & b_2 & c_6
\end{array} \right)$$
and its determinant equals $ -a_1 Res_{1,1,2}.$ The extraneous factor is 
the $1\times 1$ minor formed by the element in the fourth row, second 
column.
\par
On the other hand, because of Lemma \ref{square},  we can exhibit a 
determinantal formula for $\pm Res_{1,1,2},$ and it is given by Proposition 
\ref{small} for $t=\left[\frac{t_3}{2}\right] =0$
by the determinant of
$$
\left(
\begin{array}{ccc}
 \Delta_{(1,0,0)} & a_1 & b_1 \\
 \Delta_{(0,1,0)}& a_2 & b_2 \\
 \Delta_{(0,0,1)}& a_3 & b_3 \\
\end{array}
\right),
$$
where $\Delta_\gamma$ are  coefficients of the  Bezoutian (\ref{delta}).
Explicitly, we have 
$$\Delta_{(1,0,0)} = c_1 (a_2 b_3 -  a_3 b_2)  - c_{4} (a_1 b_3  - 
a_{3} b_{1})  + c_5 (a_1 b_2  - a_{2} b_{1}),$$
$$ \Delta_{(0,1,0)} = c_{6} (a_1 b_2  - a_{2} b_{1}) - c_{2} (a_1 b_3  - 
b_1 a_3)$$
and 
$$\Delta_{(0,0,1)}=c_{3} (a_1 b_2  - b_1 a_2).$$
This is the matrix $M_0$ corresponding to the linear transformation $\Psi_0$
which is defined in (\ref{prindis}).

\smallskip

Take now $n=4,$ and $\left(d_1,d_2,d_3,d_4\right)=\left(1,1,2,3\right).$
The critical degree is $ 3.$ 
Macaulay's classical matrix, $M_4,$  has size  $35\times 35.$ 
Because the degree of 
$\rm Res_{1,1,2,3}$ is $2+3+6+6=17,$ we know that its
extraneous factor must be a minor of size $18\times 18.$ 
By Proposition \ref{small}, we can find the smallest possible matrix
for $t= 1$ or $t=2.$ Set $t=2.$
We get the following $12\times 12$ matrix 
$$\left(
\begin{array}{cccccccccccc}
\Delta^1_{(2,0,0,0)}&\Delta^2_{(2,0,0,0)}&\Delta^3_{(2,0,0,0)}&\Delta^4_{(2,0,0,0)}&
a_1 &0&0 &0 & 0 &0 & 0 & c_1 \\
\Delta^1_{(0,2,0,0)}&\Delta^2_{(0,2,0,0)}&\Delta^3_{(0,2,0,0)}&
\Delta^4_{(0,2,0,0)}&0 &a_2 &0 &0 & b_2 &0 & 0 & c_2 \\
\Delta^1_{(0,0,2,0)}&\Delta^2_{(0,0,2,0)}&\Delta^3_{(0,0,2,0)}&
\Delta^4_{(0,0,2,0)}&0 &0 &a_3 &0 & 0 &b_3 & 0 & c_3 \\
\Delta^1_{(0,0,0,2)}&\Delta^2_{(0,0,0,2)}&\Delta^3_{(0,0,0,2)}&
\Delta^4_{(0,0,0,2)}&0 &0 &0 &a_4 & 0 &0 & b_4 &  c_4\\
\Delta^1_{(1,1,0,0)}&\Delta^2_{(1,1,0,0)}&\Delta^3_{(1,1,0,0)}&
\Delta^4_{(1,1,0,0)}&a_2 &a_1 &0 &0 & b_1 &0 & 0 & c_5 \\
\Delta^1_{(1,0,1,0)}&\Delta^2_{(1,0,1,0)}&\Delta^3_{(1,0,1,0)}&
\Delta^4_{(1,0,1,0)}&a_3 &0 &a_1 &0 &0  &b_1 & 0 & c_6 \\
\Delta^1_{(1,0,0,1)}&\Delta^2_{(1,0,0,1)}&\Delta^3_{(1,0,0,1)}&
\Delta^4_{(1,0,0,1)}&a_4 &0 &0 &a_1 &0  &0 & b_1 & c_7 \\
\Delta^1_{(0,1,1,0)}&\Delta^2_{(0,1,1,0)}&\Delta^3_{(0,1,1,0)}&
\Delta^4_{(0,1,1,0)}&0 &a_3 &a_2 &0 & b_3 &b_2 & 0 & c_8 \\
\Delta^1_{(0,1,0,1)}&\Delta^2_{(0,1,0,1)}&\Delta^3_{(0,1,0,1)}&
\Delta^4_{(0,1,0,1)}&0 &a_4 &0 &a_2 & b_4 &0 & b_2 & c_9 \\
\Delta^1_{(0,0,1,1)}&\Delta^2_{(0,0,1,1)}&\Delta^3_{(0,0,1,1)}&
\Delta^4_{(0,0,1,1)}&0 &0 &a_4 &a_3 & 0 &b_4 & b_3 & c_{10} \\
a_1&a_2&a_3&a_4 & 0&0&0&0&0&0&0&0\\
b_1&b_2&b_3&b_4 & 0&0&0&0&0&0&0&0\\
\end{array} \right)$$
where
$$
\begin{array}{ccc}
f_1 &= &a_1 X_1 + a_2 X_2 + a_3 X_3 + a_4 X_4\\
f_2 &=& b_1 X_1 + b_2 X_2 + b_3 X_3 + b_4 X_4\\
f_3 &=& c_1 X_1^2 + c_2 X_2^2 + c_3 X_3^2+ c_4 X_4^2+ c_5 X_1 X_2 + c_6 
X_1 X_3 \\
&&+ c_7 X_1 X_4 + c_8 X_2 X_3 + c_9 X_2 X_4 + c_{10} X_3 X_4,\\
\end{array}
$$
$f_4$ is a homogeneous generic polynomial of degree $3$ in
four variables, and for each $\gamma, \, |\gamma|=2,$  we write
$$\Delta_\gamma(X) = \sum_{j=1}^4 \Delta_\gamma^j X_j,$$
which has degree $1$ in the coefficients of each $f_i, i=1,\dots,4.$
The determinant of this matrix  is actually
$\pm a_1 {\rm Res}_{1,1,2,3}.$ Here, 
the extraneous factor is the minor $1\times 1$ of the matrix obtained by
taking the element in the fifth row, sixth column. 

\smallskip
In the following table, we display the minimal size  of 
the matrices $M_{t}$ and the size of classical Macaulay matrix 
 for several values of $n, \ d_1,\dots,d_n.$

$$
\boxed{
\begin{array}{cccccccc}
n && \left(d_1,\dots,d_n\right) && {\text { min size }} && 
{\text { classical }} \\
&& && && \\
2 && \left(10,70\right) && 70 && 80 \\
2 && \left(150,200\right) && 200 &&  350\\
3 && \left(1,1,2\right)&& 3 && 6\\
3 && \left(1,2,5\right)&& 14 && 28\\
3 && \left(2,2,6\right)&& 21 && 45 \\
4 && \left(1,1,2,3)\right)&& 12&&35\\
4 && \left(2,2,5,5\right) && 94 && 364 \\
4 && \left(2,3,4,5\right)&& 90 && 364\\
5 && \left(4,4,4,4,4\right) && 670 && 4845 \\
7 && \left(2,3,3,3,3,3,3\right)&& 2373 && 38760\\
10&& \left(3,3,\dots,3\right) && 175803 && 14307150\\
20&& \left(2,2,\dots,2\right) && 39875264&&131282408400\\
\end{array}
}.
$$
We give in section 4 an estimate for the ratio between these sizes. 
However, it should be noted that the number of coefficients of the 
Bezoutian that one needs to compute increases when the size of
the matrix $M_t$  decreases. We refer to \cite{fgs} and \cite{sax}
for complexity considerations  on the computation of 
Bezoutians. In particular, this computation can be well parallelized.
Also, the particular structure of the matrix and the coefficients
could be used to improve the complexity estimates; this problem is studied 
for $n=2$ and $n=3$ in \cite{czh}.
\medskip

Our approach  combines Macaulay's  original ideas  \cite{Mac}, 
expanded by Jouanolou
in \cite{jou}, with the expression for the resultant as the determinant of
a Koszul complex inspired by the work of Cayley \cite{cay}. 
We also use the work  
\cite{Ch1}, \cite{cha} of Chardin on homogeneous subresultants,
where a Macaulay 
style formula for subresultants is presented. 
In fact, we show that the proposed determinants are explicit 
non-zero minors of a bigger
matrix which corresponds to  one of the morphisms
in a Koszul resultant complex which in general has many non zero terms,
and whose determinant is $\Rd(\efen)$ (cf. Theorem \ref{rescth}).
These are the complexes considered in \cite{wz},\cite{gkz} in the equal 
degree case, built from the spectral sequence associated 
with a twisted Koszul complex at the level of sheaves.
\par
We give explicit expressions 
for the morphisms in these complexes in terms of the  Bezoutian
associated with $\efen$  for  degrees under critical degree, addressing
in this manner a problem raised by Weyman and Zelevinsky  in \cite{wz}
(cf. also \cite[13.1.C]{gkz}.
\par
In the last sections, we show  that  different classical 
formulas can be
viewed as special cases of the determinantal formulas that we present here
(cf.  \cite{gkz},\cite{wz}). In particular, we also recover  in this setting the
``affine'' Dixon formulas considered in \cite{EM} and we 
classify in particular all such determinantal formulas.

\medskip

\section{ Notations and some preliminary statements}

Let  $S_u$ denote the $A$-free module generated by the monomials in 
$A[X]$ with degree $u.$ If $u<0,$ then we set $S_u= 0.$ 
Define also the following free submodules 
$E^{t,j} \, \subseteq \, S^{t,j} \, \subseteq \, S_{t-d_j},$
for all  $j=1,\dots,n:$
\begin{equation}
S^{t,j}:= \left<X^\gamma, |\gamma|=t-d_{j},
\gamma_1<d_{1},\dots, \gamma_{j-1}< d_{j-1}\right>  
\end{equation}
\begin{equation}
E^{t,j}:= \left<X^\gamma\in S^{t,j},
\mbox{there exists } \ i\neq j: \gamma_i\geq d_{i}\right>.
\end{equation}
Note that $E^{t,n}=0,$ and $S^{t,1}=S_{t-d_1}\
\forall t\in\NN_0.$
\par
Let $j_u:S_u\rightarrow  S_u^*$ be the  isomorphism associated with
the monomial bases in $S_u$ and denote by $T_\gamma:=j_u(X^\gamma)$ the
elements in the dual basis.

\smallskip

\begin{con}
All spaces that we will consider have a monomial basis, or a dual monomial
basis. We shall suppose all these bases have a fixed order. 
This will allow us to define  matrices ``in the monomial bases'', with
no ambiguity.
\end{con}

\smallskip

Let $\psi_{1,t}$ be the $A$-linear map
$$\psi_{1,t}: S_{t_n-t}^* \rightarrow   S_{t} $$
 which sends 
\begin{equation}
\label{eq:psi1}
T_\gamma  \mapsto \Delta_{\gamma}\left(X\right),
\end{equation}
where the polynomial $\Delta_{\gamma}\left(X\right)$ is 
defined  in (\ref{delta}). Let $\Delta_t$ denote the  matrix of $\psi_{1,t}$
in the monomial bases.
\medskip
\begin{lemma}
\label{bezut}
For suitable orders of the monomial bases
in $S_t$ and $S_{t_n-t},$ we have that
$$
^{\bf t}\Delta_t=\Delta_{t_n-t}.
$$
\end{lemma}
\begin{proof}
It holds  that $\Delta\left(X,Y\right)=\Delta\left(Y,X\right)$ by the
symmetry property of Bezoutians (cf. \cite[3.11.8]{jou}).
This implies that $$\sum_{|\gamma|=t_n-t}\Delta_\gamma\left(X\right) Y^\gamma=
\sum_{|\lambda|=t}\Delta_\lambda\left(Y\right) X^\lambda=
\sum_{|\gamma|=t_n-t,\  |\lambda|=t}c_{\gamma\lambda}X^\lambda Y^\gamma,$$
with $c_{\gamma\lambda}\in A.$
It is easy to see that if
$\Delta_t=\left(c_{\gamma\lambda}\right)_{|\gamma|=t_n-t, \ |\lambda|=t}$
then 
$\Delta_{t_n-t}=\left(c_{\gamma\lambda}\right)_{|\lambda|=t, \
|\gamma|=t_n-t}.$
\end{proof}

Let us consider also the {\it Sylvester \/} linear map $ \psi_{2,t}:$
\begin{equation}
\label{prev}
\begin{array}{cccccccc}
\psi_{2,t}: & S^{t,1} & \oplus  & \cdots  & \oplus  & S^{t,n}   & \rightarrow  
& S_{t} \\ 
 & (\, g_{1} & , & \dots  & , & g_{n} \, ) & \mapsto  & 
\sum_{i=1}^{n}g_{i}f_{i},
\end{array}
\end{equation}
and denote by $D_t$ its matrix in the monomial bases. 
As usual, $\psi_{2,t_n-t}^*$
denotes the dual mapping of (\ref{prev})  in degree $t_n-t.$

Denote
\begin{equation}
\label{prindis}
\Psi_t: S_{t_n-t}^* \oplus \left(S^{t,1}\oplus 
\cdots \oplus S^{t,n} \right)  
 \rightarrow  S_{t} 
\oplus \left (S^{t_n-t,1}  \oplus   \cdots \oplus  
S^{t_n-t,n}\right)^*
\end{equation}
the $A$-morphism defined by
\begin{equation}
\label{prinmap} 
( T\, , \, g \ )  \mapsto   
\ (\psi_{1,t}(T)+\psi_{2,t}(g), \psi_{2,t_n-t}^*\left(T)\right),
\end{equation}
and call $M_t$ the matrix of $\Psi_t$ in 
the monomial bases. 

Denote also by $E_t$ the submatrix of $M_t$ whose columns are 
indexed by
the monomials in $E^{t,1}\cup \dots \cup E^{t,n-1},$ and  whose rows
are indexed by the monomials 
$X^\gamma$
in $S_t$ for which there exist two different indices $i,j$ such that
$\gamma_i \geq d_i, \gamma_j\geq d_j.$ 
With these choices it is not difficult to see that
$M_t$ and $E_t$ (when defined) are square matrices. 
\begin{remark}
Observe that $E_t$ is actually a submatrix of $D_t.$
In fact, $E_t$ is transposed of the square submatrix 
named ${\mathcal E}
\left(t\right)$ in \cite{cha}, and whose determinant is denoted by
$\Delta\left(n,t\right)$ in \cite[Th. 6]{Mac}.
\end{remark}

\begin{lemma}\label{size}
$M_t$ is a square matrix of size $\rho(t),$ where $\rho$ is
the function defined in (\ref{tamatrix}).
\end{lemma}

\begin{proof}
The assignment which sends a monomial $m$ in $S^{t,i}$ to
$x_i^{d_i} \cdot m$ injects the union of the monomial bases in
 each  $S^{t,i}$ onto the monomials of degree $t$ which are divisible
by some $x_i^{d_i}.$ 
It is easy to see that the cardinality of  the set of complementary
monomials of degree $t$ is precisely $H_d(t)$,
where $H_{d}(t)$ denotes 
the dimension of the $t$-graded
piece of the quotient  of the polynomial ring over $k$  
by the ideal generated by a regular sequence
of  homogeneous polynomials with degrees $d_1,\dots,d_n$
 (cf. \cite[3.9.2]{jou}). Moreover, using the assignment 
$\left(\gamma_1,\dots,\gamma_n\right)\mapsto
\left(d_1-1-\gamma_1,\dots,d_n-1-\gamma_n\right),$  it follows that
\begin{equation}
\label{sym}
H_d(t) = H_d(t_n-t). 
\end{equation}
\par
We can compute explicitly this Hilbert function by the following
formula (cf. \cite[\S 2]{Mac}):

\begin{equation}
\label{hilbertf}
\frac{\prod_{i=1}^{n}\left( 1-Y^{d_{i}}\right) }{\left( 1-Y\right) ^{n}}=
\sum_{t=0}^\infty H_{d}(t).Y^t.
\end{equation}

Then, $${\rm rk} \, (S^{t,1}\oplus 
\cdots \oplus S^{t,n}) = {\rm rk} \, S_t - H_d(t) .$$
Similarly, 
$${\rm rk}\left(S^{t_n-t,1}\oplus 
\cdots \oplus S^{t_n-t,n}\right)^* = {\rm rk} \, (S_{t_n-t})^*
- H_d(t_n-t).$$ 

Therefore, $M_t$ is square of
size ${\rm rk} \, S_t - H_d(t_n-t) + {\rm rk} \, S_{t_n-t}.$ Since
$i(t_n-t)= {\rm rk} \, S_{t_n-t} - H_d(t_n-t),$ the size of
$M_t$ equals
${\rm rk} \, S_t  + i(t_n-t) = \rho(t).$
\end{proof}

\begin{remark} Ordering properly 
the monomial bases, $M_t$ is the transpose of the matrix
which appears in \cite[3.11.19.7]{jou}. It has  the following structure:
\begin{equation}
\label{mjouanolou}
 \left[ 
\begin{array}{cc}
\Delta_{t} &D_{t} \\ 
 ^{\bf t} D_{t_n-t} &  0 \\
\end{array}
\right].
\end{equation}
\end{remark}
\begin{remark}
Because $\psi_{2,t}=0$ if and only if $t<\min\{d_i\},$ we have that 
$\Psi_t=\psi_{2,t}+\psi_{1,t}$ if $t>t_n-\min\{d_i\},$
and  $\Psi_t=\psi_{2,t}$ if $t>t_n.$
\end{remark}
Finally, denote ${\EE}_t$ the square submatrix of $M_t$ which has the following
structure:
\begin{equation}
\label{efactor}
 {\EE}_t=\left[ 
\begin{array}{cc}
* &E_{t} \\ 
^{\bf t} E_{t_n-t} &  0 \\
\end{array}
\right].
\end{equation}
It is clear from the definition 
that $\det({\EE}_t)=\pm\det(E_t)\det(E_{t_n-t}).$
\begin{remark}
Dualizing (\ref{prinmap}) and using lemma \ref{bezut} with a careful 
inspection at (\ref{mjouanolou}) and (\ref{efactor}), we have that  
ordering properly  their rows and
columns,
$$^{\bf t}M_t=M_{t_n-t} \ \ \mbox{and} \ \ ^{\bf t}\EE_t=\EE_{t_n-t} . $$
\end{remark}

\medskip

\section{Generalized Macaulay formulas}

We can extend the map $\psi_{2,t}$ in (\ref{prev}) to the direct
sum of all  homogeneous polynomials with degrees $t-d_1,\dots, t-d_n,$
and the map $\psi_{2,t_n-t}$ to the direct sum of all  homogeneous
polynomials with degrees $t_n-t-d_1,\dots, t_n-t-d_n,$ to get  a map
$$
\tilde{\Psi}_t: \left(S_{t_n-t}\right)^* \oplus 
\left(S_{t-d_1}\oplus \cdots \oplus S_{t-d_n} \right)  
\rightarrow  S_{t} 
\oplus \left (S_{t_n-t-d_1}  \oplus   \cdots \oplus  
S_{t_n-t-d_n}\right)^*.$$
We can thus see the matrix $M_t$ of  $\Psi_t$ in (\ref{prindis}) as a choice 
of a square submatrix of $\tilde{\Psi}_t.$ 
We will show that its determinant is a non zero minor of maximal size.

\begin{proposition}
\label{util}
Let $M'_t$ be a square matrix over $A$ of the form
\begin{equation}
\label{prima}
M'_t:=\left[ 
\begin{array}{cc}
\Delta_{t} & F_t \\ 
 ^{\bf t}F_{t_n-t}&  0 \\
\end{array}
\right]
\end{equation}
where $F_t$ has $i(t)$ columns and corresponds to a restriction of the
map
$$\begin{array}{ccc}
S_{t-d_1}\oplus\dots\oplus S_{t-d_n}&\rightarrow&S_t\\
\left(g_1,\dots,g_n\right)&\mapsto&\sum_{i=1}^n g_i \, f_i;\\
\end{array}
$$
and similarly for $F_{t_n-t}$ in degree $t_n-t.$  
Then,  $\det(M'_t)$ is a multiple 
of $\Rd(\efen)$  (probably zero).
\end{proposition}
\begin{proof}
It is enough to mimic  for the matrix $M'_t$
the  proof performed by Jouanolou in \cite[Prop. 3.11.19.10]{jou}
 to show that the determinant of the matrix $M'_t$ is an inertia 
form of the ideal $\left<\efen\right>$ (i.e. a multiple of the resultant).
We include this proof for the convenience of the reader.        
\par
Let $N:=\sum_{i=1}^n \#\{\alpha_i\in\NN^n:
|\alpha_i|=d_i\}.$
Given an algebraically closed field $k,$ and 
$a = (a_{\alpha_i})_{ |\alpha_i|= d_i,  \  i=1,\ldots,n}, $
a point in $k^N,$ 
we denote by  $f_1(a), \dots, f_n(a)\in k[X]$
the polynomials obtained from $\efen$ when the coefficients
are specialized to $a$, and similarly for the coefficients of
the Bezoutian.
Because of the irreducibility of $\Rd(\efen),$ it is enough to
show that for all $a\in k^N$ such that $f_1(a),\ldots,f_n(a)$ 
have a non trivial solution in $k^n,$ the determinant of the
specialized matrix $M'_t(a)$ is equal to $0.$
\par
Suppose that this is case, 
and let $\left(p_1,\ldots,p_n\right)$ be a non
trivial solution. Without loss of generality, we can suppose $p_1
\neq0.$
One of the rows of $M'_t(a)$ is indexed by $X_1^t.$ Replace all the
elements in that row as follows:
\begin{enumerate}
\item if the element belongs to a column  indexed by a monomial 
$X^\gamma, \ |\gamma|=t_n-t,$ then replace it with
$\Delta_\gamma (a);$ 
\item if it belongs to a column indexed by a monomial $X^\gamma\in
S_{t-d_i},$ replace it with $X^\gamma\,f_i(a).$
\end{enumerate}
It is easy to check that, the determinant of the
modified matrix is equal to $X_1^t\,\det(M'_t(a)).$
Now, we claim that under the specialization $X_i\mapsto p_i,$ 
the determinant of the modified matrix will be equal to zero if and
only if $\det(M'_t(a))=0.$
\par
In order to prove this, we will show that the following submatrix 
of size $\left(i(t_n-t)+1\right)\times \binom{n+t-1}{n-1}$ 
has rank less or equal than $i(t_n-t):$ 
$$
\left[\begin{array}{ccc}
\Delta_{\gamma_1}(a)(p)&\ldots&\Delta_{\gamma_s}(a)(p)\\
&^{\bf t}F_{t_n-t}(a)&\\
\end{array}\right].
$$
This, combined with a Laplace expansion of the determinant of the modified
matrix, gives the desired result.
\par
If the rank of the block $\left[^{\bf t}F_{t_n-t}(a)\right]$ is  less than
$i(t_n-t),$ then the claim follows straightforwardly. Suppose this is not
the case. Then the family 
$\{X^{\gamma}\,f_i(a), \ X^\gamma\in S_{t_n-t-d_i}\}$
is  a basis of the piece of degree $t_{n}-t$ of the
generated ideal $I(a):=
\langle f_1(a),\ldots,f_n(a)\rangle.$
We will show that in this case the polynomial
$ \sum_{|\gamma|=t_n-t}\Delta_{\gamma}(a)(p)X^\gamma$ belongs to $I(a)$,
which proves the claim.
\par
Because of (\ref{deltaij}) and (\ref{delta}), 
the polynomial $\left(X_1-Y_1\right)\Delta(a)(X,Y)$ lies in the ideal
$\langle f_1(a)(X)- f_1(a)(Y),\ldots,f_n(a)(X) -f_n(a)(Y)\rangle$.
Specializing $Y_i\mapsto p_i,$ we deduce that
$ (X_1-p_1) \sum_{j=0}^{t_n}(\sum_{|\gamma|=j}
\Delta_{\gamma}(a)(p)\, X^\gamma)$ is in the graded ideal $I(a).$
This, combined with the fact that $p_1\neq0,$ proves that
$ {\sum_{|\gamma|=j}
\Delta_{\gamma}(a)(p)\, X^\gamma}\in I(a)$ for all $j.$
\end{proof}
\medskip

In particular, $\Rd(\efen)$ divides $\det(M_{t}).$
We describe the extraneous factor explicitly in the following theorem,
which is the main result in this section.
Before stating it, we set
the following convention:
if the matrix ${\EE}_t$ is indexed by an empty set, 
 we define $\det\left({\EE}_t\right)=1.$

\begin{theorem} 
\label{mainth}
For any $t\geq0,$ 
 $\det\left(M_t\right) \neq 0$ and $\det(\EE_t) \neq 0.$
\par
Moreover,  we have the following formula {\it \`a la Macaulay \/}:
$$ \Rd(\efen) = \pm \frac{\det(M_t)} {\det(\EE_t)} \, .$$ 
\end{theorem}
\medskip
For the proof of  Theorem \ref{mainth}, 
we will need the following auxiliary lemma. Let
$D_t$ and $E_t$ be the matrices defined in
\S 2 before Lemma \ref{size}.

\begin{lemma}
\label{lemauxiliar}
Let $t\geq0,$ and $\Lambda$ a ring which contains $A.$ Suppose we have a
square matrix $M$ with coefficients in $\Lambda$ which has the following
structure: 
$$ M=\left[ 
\begin{array}{cc}
M_1 & D_t\\ 
 M_2 &  0 \\
\end{array}
\right],
$$
where $M_1, M_2$ are rectangular matrices. 
Then, there exists an element $m\in\Lambda$ such that
$$\det\left(M\right)=m \ .\ \det\left(E_t\right)$$
\end{lemma}
\begin{proof}
$D_t$ is square if and only if $t>t_n.$ (cf.
\cite[\S 3]{Mac}).
In this case,
$$\det(M)= \pm \det(M_2)\det(D_t);$$
because of
Macaulay's formula (cf. \cite[Th. 5]{Mac}), we
have that the right hand side equals
$$\pm\det(M_2) \det(E_t)
\Rd(\efen),$$
and the conclusion follows easily.
\par
Suppose now
$0\leq t\leq t_n.$ As in the introduction, let $i(t)$ denote the dimension
of the $k$-vector space  
of elements of degree $t$ in the ideal generated by a
regular sequence  of $n$ polynomials with degrees $d_1,\dots,d_n.$
 Then $D_t$ has  $i(t) + H_d\left(t\right)$ rows and $i(t)$  columns,
and  there is a bijection between the family  ${\F}$ of $H_d\left(t\right)$ 
monomials of degree $t,$ and the maximal minors $m_{\F}$ of $D_t.$ 
Namely, $m_{\F}$ is the determinant of the square submatrix made
by avoiding all rows indexed by monomials in ${\F}.$
\par
It is not hard to check that $m_{\F}$ is the determinant
$\phi_{\F}^*$ which is used in  
\cite{cha}, for computing the {\it subresultant} associated with the
family $\{X^\gamma\}_{\gamma\in{\F}}.$
\par
Now, using the generalized Macaulay's formula for the subresultant 
(cf. \cite{cha}), we have that
$$m_{\F}= \pm \det(E_t). \Delta^t_{\F},$$
where $\Delta^t_{\F}$ is the subresultant associated with the family 
${\F}.$ It is a polynomial in $A$ which vanishes under a
specialization of the coefficients $ f_1(a),\dots,  f_n(a)$ if and only
if the family $\{X^\gamma\}_{\gamma\in{\F}}$ fails to be  a basis of the 
$t-$ graded piece of the quotient 
$k[X_1,\dots,X_n]/  \langle f_1(a), \dots, f_n(a) \rangle$ (cf. \cite{Ch1}).
 \par
Let $m^c_{\F}$ be
the complementary minor of $m_{\F}$ in $M$ (i.e. the determinant of the
square submatrix of $M$  which is made by deleting all rows and columns that
appear in $m_{\F}$).
By the Laplace expansion of the determinant, we have that
$$
\det \left( M \right) =\sum_{\F} s_{\F} \cdot m_{\F}  \cdot m_{\F}^{c}
= \det\left(E_t\right) \left( \sum_{\F} s_{\F} \cdot m^c_{\F}  \cdot \Delta^t_{\F} 
\right)$$
with $ s_{\F}= \pm 1.$
Setting $m= \sum_{\F} s_{\F}  \cdot m^c_{\F}  \cdot \Delta^t_{\F} 
\in\Lambda,$ we have the desired result. 
\end{proof}

\medskip
We now give the proof of Theorem \ref{mainth}.

\begin{proof}
 In \cite{Mac}  it is  shown that $\det\left(E_t\right)
\neq 0, \ \forall t \geq 0.$ This implies that $\det\left(\EE_t\right)\neq0.$
In order to prove that  $\det\left(\EE_t\right)=
\det(E_t)\det(E_{t_n-t})$ divides 
$\det(M_t),$ we  use the following trick:
consider the ring $B:={\ZZ}\left[ b_{\alpha _{i}}
\right]_{|\alpha _{i}|= d_{i},i=1,\dots ,n},$ where $b_{\alpha_i}$ are
new variables, and the polynomials
$$ f_{b,i}:=\sum_{|\alpha _{i}|=d_{i}}b_{\alpha_{i}}X^{\alpha _{i}} \quad
\in B\left[ X_{1},\dots ,X_{n}\right]. $$
Let $D_t^b$ the matrix of the linear transformation 
$\psi^b_{2,t}$  determined by the formula  (\ref{prev}) 
but associated with the sequence
$f_{b,1},\dots,f_{b,n}$ instead of $\efen.$
Set $\Lambda:= \ZZ\left[a_{\alpha_i},b_{\alpha_i}\right]$, and
consider the matrix $M(a,b)$ with coefficients in $\Lambda$ given by
$$ M(a,b)=\left[ 
\begin{array}{cc}
\Delta_t & D_t\\ 
 ^{\bf t} D_{t_n-t}^b &  0 \\
\end{array}
\right].
$$
It is easy to see that $M(a,a)= M_t$, and because
of  Lemma \ref{lemauxiliar}, 
we have that $\det\left(E_t\right)$ divides $\det\left(M(a,b)\right)$ in
$\Lambda.$ Transposing $M(a,b)$ and using a symmetry argument, 
 again by  the same lemma, we can
conclude that $\det\left(E^b_{t_n-t}\right)$ divides $\det\left(M(a,b)\right)$ in
$\Lambda,$ where $E^b_{t_n-t}$ has the obvious meaning.
\par
The ring $\Lambda$ is a factorial domain and  $\det(E_t)$ and 
$\det(E_{t_n-t}^b)$
have no common factors in $\Lambda$  because they depend on
different sets of variables. So, we have 
$$\det\left(M(a,b)\right)=  p(a,b) \, \det(E_t) \, 
\det(E_{t_n-t}^b)$$
for some $p \in \Lambda.$  Now, specialize 
$b_{\alpha_i}\mapsto a_{\alpha_i}.$ 
The fact that $\det\left(M_t\right)$ is a multiple of the resultant 
has been proved in Proposition \ref{util} (see also \cite[Prop. 3.11.19.21]{jou})
for $0\leq t\leq t_n,$ and in \cite{Mac} for
$t>t_n.$  On the other side, since
$\Rd(\efen)$ is
irreducible and depends on all the coefficients of $\efen$
while $\det(E_t)$ and $ \det(E_{t_n-t})$ do not depend on the coefficients of $f_n,$
we conclude that $\Rd(\efen)$ divides $p(a,a).$ Moreover, the following lemma
shows that they have the same degree. Then, their ratio is a rational 
number $\lambda.$
We can see that $\lambda=\pm 1,$  considering the specialized family
$X_1^{d_1}, \dots, X_n^{d_n}.$
\end{proof}

\begin{lemma}
For each $i=1,\dots,n$ the degree $\deg_{(a_{\alpha_i})}\left(M_t\right)$
of $M_t$ in the coefficients of $f_i$ equals
$$
\begin{array}{ccc}
 \deg_{(a_{\alpha_i})}\left(\Rd(\efen)\right)+ \deg_{(a_{\alpha_i})}\left
(E_t\right)+\deg_{(a_{\alpha_i})}\left(E_{t_n-t}\right) &=&\\
 d_1\dots d_{i-1}.\, d_{i+1}\dots d_n +\deg_{(a_{\alpha_i})}\left
(E_t\right)+\deg_{(a_{\alpha_i})}\left(E_{t_n-t}\right)
 &&
\end{array}
$$
\end{lemma}
\begin{proof}
Set $J_u(i):= \{X^\gamma\in S_u,
 \gamma_i\geq d_i, \, \gamma_j<d_j \, \forall j\neq i\}, \ u = t,t_n-t.$
{}From the definitions of $\psi_{2,t}$ and $E_t,$ it is easy to check that,
if $\delta_t$ is a maximal minor of $D_t,$
$$ \deg_{(a_{\alpha_i})}\left(\delta_t\right)
- \deg_{(a_{\alpha_i})}\left(E_{t}\right)= \# J_t(i).$$

Using Laplace expansion, it is easy to see that $\det\left(
M_t\right)$ may be expanded as follows
$$\det\left(M_t\right) = \sum_{\delta_t,\delta_{t_n-t}} s_\delta \cdot
m_\delta \cdot  \delta_t  \cdot  \delta_{t_n-t} $$ 
where $s_\delta = \pm1$, $\delta_{t_n-t}$ is a maximal minor of $^{\bf t} D_{t_n-t}$
and $m_\delta$ is a minor of size $H_d(t)$ in $\Delta_{t}.$

As each entry of $\Delta_{t}$ has degree $1$ in the coefficients 
of $f_i,$ the lemma will be proved if we show that 
\begin{equation}
\label{sumacard}
\# J_t(i) + \# J_{t_n-t}(i) + H_d(t) = d_1\dots d_{i-1} . \, d_{i+1}\dots d_n.
\end{equation}
Now, as already observed in the proof of Lemma \ref{size},
 $H_d(t)$ can be computed as the cardinality of the following set:
\begin{equation}
\label{hdt} 
H_{d,t}:=\{X^\gamma\in S_t, \,  \gamma_j < d_j \, \forall j\},
\end{equation}
and $d_1\dots d_{i-1}. \, d_{i+1}\dots d_n$ is 
the cardinality of
$$\Gamma_{i}:=\{X_1^{\gamma_1}\dots X_{i-1}^{\gamma_{i-1}} X_{i+1}^{\gamma_{i+1}}\dots
X_n^{\gamma_n}\, , \, \gamma_j < d_j \, \forall j\}.$$
In order to prove (\ref{sumacard}) it is enough to exhibit a bijection
between $\Gamma_{i}$ and the disjoint union $J_t(i)\bigcup J_{t_n-t}(i)\bigcup
H_{d,t}.$ This is actually a disjoint union for  all $t,$ unless
$t_n-t=t.$ But what follows shows that the bijection is well defined 
even in this case.

Let $X^{\widehat \gamma}\in \Gamma_{i}, \, \widehat \gamma=
\left(\gamma_1,\dots,\gamma_{i-1},\gamma_{i+1},\dots
,\gamma_n \right)$ with $\gamma_j<d_j \, \forall j\neq i.$
If $|\widehat \gamma|\leq t,$ then there exists a unique $\gamma_i$ such that
$\gamma:= \left(\gamma_1,\dots,\gamma_n\right)\in\NN_0^n$ verifies
$|\gamma|=t.$ If $\gamma_i<d_i,$ then we send $X^{\widehat \gamma}$ to
$X^\gamma \in H_{d,t}.$ Otherwise, we send it to $X^\gamma \in J_{t}(i).$

If $|\widehat \gamma|> t,$ let $\widehat\gamma^*$ denote the
multiindex
$$\left(d_1-1-\gamma_1,\dots,d_{i-1}-1-\gamma_{i-1},
d_{i+1}-1-\gamma_{i+1},\dots, d_n-1-\gamma_n\right).$$
Then, $|\widehat\gamma^*|< t_n-t,$ and there exists 
a unique $\gamma_i$  such that the multiindex $\gamma$ defined
by
$$
\left(d_1-1-\gamma_1,\dots,d_{i-1}-1-\gamma_{i-1},\gamma_i,
d_{i+1}-1-\gamma_{i+1},\dots, d_n-1-\gamma_n\right)$$
has degree $t_n-t.$ We can send $X^{\widehat\gamma}$ to
$X^\gamma \in J_{t_n-t}(i)$ provided that $\gamma_i\geq d_i.$
Suppose this last statement does not happen, this implies that
the monomial with exponent
$$\gamma^*:=\left(\gamma_1,\dots,\gamma_{i-1},d_i-1-\gamma_i,\gamma_{i+1},\dots,
d_n\right)$$
has degree $t$ contradicting the fact that $|\widehat \gamma|> t.$

With these rules, it is straightforward to check that we obtain the desired bijection. 
\end{proof}

\medskip
Changing the order of the sequence $(\efen)$, and 
applying Theorem (\ref{mainth}), we deduce that

\begin{corollary}
\label{remate}
$\Rd(\efen) = \gcd \{\mbox{maximal minors of } \, \, \tilde\Psi_t\}.$
\end{corollary}
\medskip

\section{Estimating the size of $M_t$}
We have, for each integer $t\geq0,$ a matrix $M_t$ of size $\rho(t),$ where
$\rho$ was defined in (\ref{tamatrix}), whose determinant is
a nontrivial multiple of the resultant, and such that, moreover, its extraneous factor
is a minor of it. We want to know which is the smallest matrix we can
have.

We can write $\rho$ as
$$\rho(t) = \binom{n+t-1}{n-1} + \binom{n + t_n -t -1}{n-1} - H_d(t_n-t).$$
It is straightforward to check that 
$\binom{n+t-1}{n-1} + \binom{n + t_n -t -1}{n-1}$ is the restriction to the integers of 
a polynomial $\phi(t)$ in a real variable $t, $ symmetric with respect to $\frac{t_n}{2}$ 
(i.e. $\phi(\frac{t_n}{2}+t)=\phi(\frac{t_n}{2}-t)$ for all $t$). 
Moreover, $\phi$
reaches its minimum over $\left[0,t_n\right]$ at $t =\frac{t_n}{2}.$
Since
\begin{equation}
\label{rho}
\rho(t) =  \phi(t) -H_d(t) = \phi(t_n-t) - H_d(t_n-t) = \rho(t_n-t),
\end{equation}
in order to  study the behaviour of $\rho$
we need to understand how $H_d(t)$ varies with $t$.
We denote as usual the integer part of a real number $x$ by the symbol
$\left[x\right].$

\begin{proposition}
\label{nondec}
$H_d(t)$ is non decreasing on (the integer points of )
the interval $\left[0,\left[\frac{t_n}{2}\right]\right].$ 
\end{proposition}

\begin{proof}
We will prove this result by induction on $n$. The case
$n=1$ is obvious since $t_1=d-1$ and $H_d(t) = 1$ for any $t =0,\dots,d-1.$
Suppose then that the statement holds for $n$ variables and set
$$\widehat d:=\left(d_1,\dots,d_{n+1}\right)\in\NN_0^{n+1},$$
$$d:=(d_1,\dots,d_n).$$
Let $t < t+1 \leq \left[\frac{t_{n+1}}{2}\right].$ We want to
see that $\varphi(t):= H_{\widehat d}(t+1) - H_{\widehat d}(t)$
is non negative. 
Recall from (\ref{hdt}) that, for every $t\in\NN_0$, 
$H_{\widehat d}(t) $  equals the cardinality of the
set  
$$\{\gamma\in\NN_0^{n+1}: \ |\gamma|=t,
\, 0\leq \gamma_i\leq d_i-1, \ i=1,\dots,n+1 \}. $$ 
Then, it can also be computed as
$$\sum_{j=0}^{d_{n+1}-1}\# \{\widehat\gamma\in\NN_0^n: \ 
|\widehat\gamma|=t-j,
\, 0\leq \widehat\gamma_i\leq d_i-1, \ i=1,\dots,n 
\},$$
which gives the equality
 $H_{\widehat d}(t)= \sum_{j=0}^{d_{n+1}-1} H_d(t-j).$ 
It follows that $\varphi(t) = H_d(t+1) - H_d(t+1-d_{n+1}).$
\par
If $t+1 \leq \left[\frac{t_{n}}{2}\right],$ we deduce that $\varphi(t)
\geq 0$ by inductive hypothesis. Suppose then that $t+1$ is in the
range
$\left[\frac{t_{n}}{2}\right] < t  +1 \leq \left[\frac{t_{n+1}}{2}\right].$
As $H_d(t+1)=  H_d(t_n -t-1),$  it
is enough to show that $t_n -t -1 \geq t+1 -d_{n+1}$ and
 $t_n -t-1 \leq \left[\frac{t_{n}}{2}\right],$ which can be
easily checked, and the result follows again by inductive
hypothesis.
\end{proof}

\begin{corollary}
\label{min}
 The size
$\rho\left(t\right)$ of the matrix $M_t$ is minimal over $\NN_0$ 
when $t=\left[\frac{t_n}{2}\right].$
\end{corollary}

\begin{proof} By (\ref{rho}), $\rho$ has a maximum at 
$\left[\frac{t_{n}}{2}\right]$ over $[0,t_n]$ because $\phi$
has a maximum and $H_d$ has a minimum. 
If $t>t_n,$ we have that
$\rho\left(t\right)=\binom{n+t-1}{n-1}.$
For $t$ in this range, it is easy to check that 
$\rho\left(t_n\right)=\binom{n+t_n-1}{n-1}-1<
\rho\left(t\right).$ Then, $\rho(t) > \rho(t_n)
\geq \rho\left(\left[\frac{t_n}{2}\right]\right).
$
\end{proof}

\begin{remark}
Note that when
$t_n$ is odd, $\rho(\left[\frac{t_{n}}{2}\right]) =
\rho(\left[\frac{t_{n}}{2}\right]+1)$, and then the size
of $M_{t}$ is also minimal for $t=\left[\frac{t_{n}}{2}\right]+1$
in this case.
\end{remark}

Denote $p:= \frac{\sum_{i=1}^n d_i}{n}$ the average value of the degrees,
and set $q:=\frac{p+1}{2 p}.$  Note that except in the linear case when
all $d_i=1,$ it holds that $p>1$ and $q<1.$

\begin{proposition}
\label{rate}
Assume $p >1.$
The ratio between the size of the smallest matrix $M_t$ and the
classical Macaulay matrix $M_{t_n+1}$ can be bounded by
$$
\frac{\rho\left(\left[t_n/2\right]\right)}{\rho\left(t_n+1\right)}
\leq 2 \, q^{n-1}.$$
In particular,  it tends to zero exponentially in $n$ when the number of 
variables tends to infinity and $p$ remains bigger that a constant
$c > 1.$.
\end{proposition}

\begin{proof}
 When $t_n$ is even, $t_n - [t_n/2] = [t_n/2]$ and when
$t_n$ is odd, $t_n -[t_n/2] = [t_n/2] +1.$ In both cases,
$$
\frac{\rho\left(\left[t_n/2\right]\right)}{\rho\left(t_n+1\right)}
\leq
\frac{2\, \binom{n+\left[t_n/2\right]}{n-1} }{
\binom{n+t_n} {n-1}}=
2\, \frac{\left(\left[t_n/2\right]+n\right)\dots
\left(\left[t_n/2\right]+2\right)}{\left(t_n+n\right)\dots
\left(t_n+2\right)}=
$$
$$
= 2\, \left(\frac{\left[t_n/2\right]+n}{t_n+n}\right)
\left(\frac{\left[t_n/2\right]+n-1}{
t_n+n-1}\right)\dots
\left(\frac{\left[t_n/2\right]+2}{
t_n+2}\right)\leq $$
$$\leq
2\, \left(\frac{\left[t_n/2\right]+n}{
t_n+n}\right)^{n-1}.$$
Since $t_n = np -n,$ we deduce that
$$\frac{\left[t_n/2\right]+n}{ 
t_n+n }\, \leq \, \frac{ \frac{np}{2} + \frac{n}{2}}{np} =
\frac{1}{2} + \frac{1}{2p} = q,
$$
as wanted.
\end{proof}

\medskip

\section{Resultant complexes}

In this section we consider Weyman's complexes (cf. \cite{wz}, \cite{gkz})
and we make explicit  the morphisms in these complexes, 
which lead to polynomial expressions for the resultant via determinantal
formulas in the cases described in Lemma \ref{square}.
\par
We will consider a complex which is a ``coupling'' of the
Koszul complex $ {\bf K}^\bullet (t;\efen)$ 
associated with $\efen$ in degree $t$ and the dual of the 
Koszul complex ${\bf K}^\bullet (t_n-t,\efen)^*$ 
associated with $\efen$ in degree $t_n -t.$
This complex  arises from the spectral 
sequence derived from the Koszul complex  of
sheaves  on $\CP^{n-1}$ associated with $\efen$
twisted by ${\mathcal O_{\CP^{n-1}}}(t).$ Here, ${\mathcal 
O_{\CP^{n-1}}}(t)$ denotes as usual the $t$-twist of the sheaf
of regular functions over the $(n-1)$-projective space $\CP^{n-1}$
(see for instance \cite[p. 34]{gkz}). Its space of
global sections can be identified with the space of homogeneous
polynomials in $n$ variables of degree $t$.
We make explicit  in terms of the Bezoutian
the map $\partial_0$ (see (\ref{prinmap}) below)
produced by cohomology obstructions. In fact, the non-trivial
contribution is given in terms of the mapping
$\psi_{1,t}$  defined in (\ref{eq:psi1}).
\par
Precisely, let ${\bf K}^\bullet (t;\efen)$ denote the complex
\begin{equation}
\label{complex1}
 \{0 \longrightarrow \, K(t)^{-n} \,  \mathop{\longrightarrow}^{\delta_{-(n-1)}}
\, \dots \, \mathop{\longrightarrow}^{\delta_{-1}} \, K(t)^{-1}\, 
\mathop{\longrightarrow}^{\delta_{0}} \, K(t)^{0}\,
\} ,
\end{equation}
where
$$K(t)^{-j} = \mathop{\oplus}_{i_1<\dots <i_j} S_{t - d_{i_1} - \dots - d_{i_j}}$$
and $\delta_{-j}$ are the standard Koszul morphisms.

Similarly,  let ${\bf K}^\bullet (t_n-t;\efen)^*$ denote the complex
\begin{equation}
\label{complex2}
 \{ 
K(t_n-t)^{0} \,  \mathop{\longrightarrow}^{\delta^*_{0}}
 \, K(t_n-t)^{1}\,
  \mathop{\longrightarrow}^{\delta^*_{1}} 
\, \dots \, \mathop{\longrightarrow}^{\delta^*_{n}} \, K(t_n-t)^{n}\,
\} ,
\end{equation}
where
$$K(t_n-t)^{j} = \mathop{\oplus}_{i_1<\dots <i_{j}} S^*_{t_n-t - d_{i_1} - \dots - d_{i_{j}}}$$
and  $\delta^*_j$ are the duals of the standard Koszul morphisms.
 Note that in fact $K(t_n-t)^n =0 $ for any $t \geq 0.$
\par
Now, define ${\bf C}^\bullet (t;\efen) $ to be the following coupled complex
\begin{equation}
\label{complex}
 \{ 0 \longrightarrow \,
C^{-n} \,  \mathop{\longrightarrow}^{\partial_{-(n-1)}}
\, \dots \, \mathop{\longrightarrow}^{\partial_{-1}} \, C^{-1}\, 
\mathop{\longrightarrow}^{\partial_{0}} \, C^{0}\,
  \mathop{\longrightarrow}^{\partial_{1}} 
\, \dots \, \mathop{\longrightarrow}^{\partial_{n-1}} \, C^{n-1} \longrightarrow 0 \,
\} ,
\end{equation}

\smallskip

\noindent where

\begin{equation}
\begin{array}{lcll}
\label{def:complex}
C^{-j} & = & K(t)^{-j}, & \, j=2, \dots, n \\
C^j & = & K(t_n-t)^{j+1}, & \, j=1,\dots,n-1\\
C^{-1} & = &  K(t_n-t)^0 \oplus  K(t)^{-1}&{}\\
C^{0} & = & K(t)^0 \oplus K(t_n-t)^1 & {}
\end{array}
\end{equation}

\smallskip

\noindent and the morphisms  are defined by

\begin{equation}
\begin{array}{lcll}
\label{pmap}
\partial_{-j} & = & \delta_{-j}, & \, j=2, \dots, n-1 \\
\partial_j & = & \delta^*_{j}, & \, j=2,\dots,n-1\\
\partial_{-1} & = & 0 \oplus   \delta_{-1} & {}\\
\partial_{0} & = &  ( \psi_{1,t}+\delta_0) \oplus \delta^*_0 & {}\\
\partial_1 & = & 0 +\delta^*_1 & {}
\end{array}
\end{equation}
More explicitly, $\partial_0( T, (g_1,\dots,g_n) )= 
( \psi_{1,t}(T) +\delta_0( g_1,\dots,g_n),  \delta^*_0 (T))$
and $\partial_1(h,( T_1,\dots,T_n)) = \delta^*_1( T_1,\dots,T_n).$
Observe that $\partial_0$ is precisely the mapping
 we called $\tilde\Psi_t$ in the previous section.

\smallskip
As in the proof of Proposition \ref{util},
given an algebraically closed field $k,$ and 
$a = (a_{\alpha_i})_{ |\alpha_i|= d_i,  \  i=1,\ldots,n}, $
a point in $k^N,$ 
we denote  by  $f_1(a), \dots,$ $f_n(a)$
the polynomials $\in k[X]$ obtained from $\efen$ when the coefficients
are specialized to $a$.
For  any particular choice of coefficients  
in (\ref{complex}) we get a complex of $k$-vector spaces.
We will denote the specialized
modules and morphisms by $K(t)^{1}(a), \delta_0(a),$ etc.
Let $D$ denote the determinant (cf. \cite[Appendix A]{gkz}, 
\cite{De}) 
of the complex of $A$-modules (\ref{complex}) 
with respect to the monomial bases
of the $A$-modules $C^\ell$. 
This is an element in the field of fractions of $A$.  
 
\smallskip

We now state the main result in this section. 

\begin{theorem}
\label{rescth}
The complex (\ref{complex}) is generically exact, and for
each specialization of the coefficients it is exact if 
and only if the resultant does not vanish. For any positive
integer $t$  we have that 
\begin{equation}
\label{det}
D =  \Rd(\efen),
\end{equation}
and moreover,  $D$  equals the
greatest common divisor of all 
maximal minors of a matrix representing the
$A$-module map  $\partial_0$. 
\end{theorem}

\begin{proof}
For $t>t_n,$ we get the Koszul complex in degree $t,$ and so the specialized 
complex at a point $a \in k^N$ is exact if 
and only if $f_1(a),\dots, f_n(a)$ is a regular sequence, i.e. if and only if 
the resultant
does not vanish.  The fact that the determinant of this complex equals the
resultant goes back to  ideas of Cayley; for a proof see \cite{De}, \cite{gkz} or \cite{Ch3}.
\par
Suppose $0\leq t\leq t_n.$  Since $\delta_0 \circ \delta_{-1} = \delta^*_1 \circ
\delta^*_0 =0,$ it is easy to see that (\ref{complex}) is a complex.
\par
Set 
$$U:= \{ a=(a_{\alpha_i}) \in k^N, i=1,\dots,n ,  |\alpha_i|=d_i : \det(M_t(a)) \not= 0 \}.$$
Note that the open set $U$ is non void because the vector of coefficients
of $\{X_1^{d_1}, \dots, X_n^{d_n}\}$ lies in $U,$ since in this case
$\det M_t= \pm 1.$ For any choice of homogeneous
polynomials $f_1(a),\dots, f_n(a) \in k[X]$ with respective degrees $d_1,\dots,d_n$ and
coefficients $a$ in $U$, the resultant does not
vanish by Theorem \ref{mainth} and then the specialized Koszul complexes
 in (\ref{complex1})  and (\ref{complex2}) are exact.   
\par
Then, the dimension
$\dim \im(\delta_{0}(a))$ of the image of $\delta_{0}(a)$
equals $i(t) = \dim < f_1(a),\dots, f_n(a) >_t.$ 
Similarly, $\dim (\ker(\delta^*_0(a)) =  i(t_n-t).$
 Therefore, 
\begin{eqnarray*}
\dim \ker(\partial_0(a)) \geq \dim \im(\partial_{- 1}(a))=
\dim \im (\delta_{-1}(a)) = \\
=\dim \ker(\delta_0(a)) =  \dim K(t)^{-1}(a) - 
i(t) .
\end{eqnarray*}
On the other side,
the fact that $M_t(a)$ is non singular of size $\rho(t)$ implies that 
\begin{eqnarray*}
\dim \ker(\partial_0(a)) \leq \dim C^{-1}(a) - \rho(t) = \\ = \dim K(t)^{-1}(a) +
\dim K(t_n-t)^0(a) - \rho(t) = \\ = \im K(t)^{-1} (a)+ \dim S_{t_n-t} (a)-\rho(t)= \\
= \dim K(t)^{-1}(a)- i(t).
\end{eqnarray*} 
\par Therefore,  $ \dim \im(\partial_{-1}(a)) = \dim \ker(\partial_0(a))$
and the complex is exact at  level $-1.$ 
\par In a similar way, we can check
that the complex is exact at level $0$, and so the full specialized complex
(\ref{complex}) is exact when the coefficients $a$ lie in  $U.$
\par
In order to compute the determinant of the complex in
this case, we can make suitable choices of monomial
subsets in each term of the complex starting from the index sets
that define $M_t(a)$ to the left and to the right.
Then,
$$D(a) = \frac{\det M_t(a)} { p_1(a) \cdot p_2(a)},$$
where $p_1(a)$ (resp. $p_2(a)$) is a quotient of
product of minors of the morphisms on the left (resp. on the right).
\par
Taking into account (\ref{complex1}) and (\ref{complex2}), it follows from
\cite{cha} that  
$$p_1 (a)= \det (E_t(a)), \, p_2(a)= \det(E_{t_n-t}(a)),$$ 
and so by Theorem \ref{mainth} we have
\begin{eqnarray*}
D(a) = \Rd(\efen)(a) \frac{\det({\EE}_t)(a)}{\det(E_t(a)) \det(E_{t_n-t}(a)) }= \\
 = \Rd(\efen)(a)
\end{eqnarray*}
for all families of homogeneous polynomials with coefficients $a$
in the dense open set $U$, and since $D$ and the resultant
are rational functions,  this implies
(\ref{det}), as wanted.  Moreover,  it  follows that the complex is exact 
if and only if the resultant does not vanish. 
\par
The fact that $\Rd(\efen)$ is the greatest common divisor of all maximal 
minors of the matrix representing $\partial_0$ has been proved in Corollary 
\ref{remate}. 
\end{proof}

\medskip
We remark that from the statement of Theorem \ref{rescth} plus a close
look at the map at level $0$, it is not hard to deduce that for a given
specialization of $\efen$ in $k$ with non vanishing resultant, the
specialized polynomials $\Delta_{\gamma}(a), |\gamma| = t_{n}-t$ 
generate the quotient of the polynomial ring $k[X]$ by the ideal
$I(a) =\langle f_{1}(a), \dots, f_{n}(a) \rangle$ in degree $t$.
We can instead use the known dualizing properties of the
Bezoutian in case the polynomials define a regular sequence, to provide
an alternative proof of Theorem \ref{rescth}. This is a consequence of
Proposition  \ref{duality} below. We refer to  
\cite{jou},[\cite{kunz},
Appendix F], \cite{ss} and \cite{ts} for the relation 
between the Bezoutian and the residue (i.e. an associated trace)
and we simply recall the properties that we will use.

\smallskip

Assume $\Rd\left(f_1(a),\ldots,f_n(a)\right)$ is  different from zero.
This implies 
that $f_1(a),\ldots,f_n(a)$ is a regular sequence and 
their zero locus consists of the single point ${\bf 0}\in k^n.$
Then, there exists a dualizing $k$-linear operator 
$$R_{0} : k[Y]/ \langle f_{1}(a)(Y), \dots,
f_{n}(a)(Y) \rangle \longrightarrow k,$$ 
called the {\it residue or trace operator\/}, which verifies
\begin{enumerate}
\item $h(X) = R_0\left(h(Y)\,\Delta(a)(X,Y)\right)$ in the
quotient ring $k[X] / I(a).$
\item If $h$ is homogeneous of degree $t$ with $t\neq t_n,$ 
$R_0(h)=0$ 
\end{enumerate}
Then, for every polynomial $h(X)\in k[X]$ of degree $t,$ it
holds that
\begin{equation}
\label{dual}
h(X)
=\sum_{|\gamma| 
=t_n-t}R_0\left(h(Y)\,Y^{\gamma}\right)\Delta_{\gamma}(a)(X) \quad
{\text { mod } } I(a),
\end{equation}
where $\Delta(a) (X,Y) = \sum_{|\gamma| = t_{n}-t} 
\Delta_{\gamma}(a)(X) Y^\gamma$ as in (\ref{delta}).
As a consequence, 
{\it the family  $\{\Delta_{\gamma}(a)(X)\}_{|\gamma|=t_n-t},$ (resp.
$|\gamma|=t$) generates the graded piece of the quotient in degree
$t$ (resp. $t_n-t$)}. Moreover, it is easy to verify that for any 
choice of polynomials $p_i(X,Y),\,q_i(X,Y)\in
k[X,Y], \ i=1,\ldots,n,$ the polynomial $\tilde{\Delta}_{a}(X,Y)$ 
defined by
\begin{equation}
\label{deltamod}
\tilde{\Delta}_{a}(X,Y):=\Delta(a)(X,Y) +\sum_{i=1}^n{
p_i(X,Y)\,f_i(a)(X) + q_i(X,Y)\,f_i(a)(Y)}.
\end{equation}
has the same dualizing properties as $\Delta(a)(X,Y)$.
\par
 
We are ready to prove a kind of ``converse'' to Proposition \ref{util}.

\begin{proposition}
\label{duality}
If $\Rd\left(f_1(a),\ldots,f_n(a)\right)\neq 0,$ it is possible
to extract  a square submatrix $M'_t$ of $\tilde{\Psi}_t$  
as in (\ref{prima}) such that $\det\left(M'_t(a)\right)\neq0.$ 
\end{proposition}
\begin{proof}
Since
$ f_1(a)(X),\ldots,f_n(a)(X)$ is a regular sequence in $k[X],$
the dimensions of the graded pieces of the quotient $k[X] / 
I(a)$ in degrees $t$ and $t_n-t$ are $i(t)$ and $i(t_n-t)$ respectively.
\par 
We can then choose blocks $F_t$ and $F_{t_n-t}$ as in (\ref{prima}) such that
$F_t(a)$ and $F_{t_n-t}(a)$ have maximal rank.
Suppose without loss of generality that the blocks $F_t$ and
$F_{t_n-t}$ have respectively the form
$\left[\begin{array}{c}
Q_t\\
R_t\end{array}\right]$ and $\left[\begin{array}{c}
Q_{t_n-t}\\R_{t_n-t} \end{array}\right]$,
where $Q_t(a)$ and $Q_{t_n-t}(a)$ are square invertible matrices of maximal 
size. 
We are going to prove that, with this choice, the matrix $M'_t(a)$ is 
invertible.
\par
Our specialized matrix will look as follows:
$$
M'_t(a) = \left[
\begin{array}{cc}
\Delta_t(a) &  {\begin{array}{c}Q_t(a)\\
R_t(a)\end{array}} \\
^{\bf t}Q_{t_n-t}(a) \ ^{\bf t} R_{t_n-t}(a)&0 
\end{array}
\right].
$$
Applying linear operations in the rows and columns of $M'_t(a),$ it can
be transformed into:
$$\left[\begin{array}{ccc}
0&0&Q_t(a)\\
0&\tilde{\Delta}_{t,a}&R_t(a)\\
^{\bf t}Q_{t_n-t}(a) & ^{\bf t} R_{t_n-t}(a)&0 
\end{array}
\right],
$$
where the block $[\tilde{\Delta}_{t,a}]$ is square and of size $H_d(t).$
\par
But it is easy to check that this $\tilde{\Delta}_{t,a}$ corresponds to the
components in degree $t$ of another Bezoutian 
$\tilde{\Delta}_{a}(X,Y)$ 
(in the sense of (\ref{deltamod})). 
This is due to the fact that each of the linear operations performed on 
$M'_t(a),$ when applied to the
block $\Delta_{t,a},$  can be read as  a polynomial combination of 
$f_i(a)(X)$ and $f_i(a)(Y)$ applied to the bezoutian $\Delta(a)(X,Y).$ 
\par
Using the fact that the polynomials $\tilde{\Delta}_{\gamma,a}(X)$
read in the columns of $\tilde{\Delta}_{t,a}$ 
generate the quotient in degree $t_n-t$ and they are as many as
its dimension, we deduce that they are a basis and so
$$
\det \left(
\begin{array}{cc}
0&\tilde{\Delta}_t \\
^{\bf t}Q_{t_n-t}(a) & ^{\bf t} R_{t_n-t}(a) 
\end{array}
\right)\neq0,
$$
which completes the proof of the claim.
\end{proof}

We could then avoid the consideration of the open set $U$ in the
proof of Theorem \ref{rescth}, and use Proposition \ref{duality}
to show directly that the complex is exact outside the zero locus
of the resultant. In fact, this is not surprising since 
for all specializations such that the resultant is non zero,
the residue operator defines a natural duality between the $t$-graded 
piece of the the quotient of the ring of polynomials
with coefficients in $k$ by the ideal $I(a)$
and the $t_n-t$ graded piece of the quotient,
and  we can read dual residue bases in the Bezoutian. 

\bigskip

We characterize now those data $n, d_1,\dots,d_n$ for
which we get a determinantal formula.

\begin{lemma}
\label{square}
Suppose $d_1 \leq d_2 \leq \dots \leq d_n.$
The determinant of the resultant complex provides a determinantal 
formula for the resultant $\Rd(\efen)$ if and only if
the following  inequality is verified
\begin{equation}
\label{squareq}
 d_3 + \dots + d_n -n < d_1 + d_2 - 1.
\end{equation}
Moreover, when (\ref{squareq}) holds, there exists a determinantal
formula given by the resultant complex for each $t$ such that
\begin{equation}
\label{squaret}
d_3 + \dots + d_n -n< t <d_1 + d_2 .
\end{equation}
\end{lemma} 

\begin{remark}
When all $d_i$ have a common value $d$, (\ref{squareq}) reads
$$(n-2) d < 2 d + n - 1,$$
which is true for any $d$ for $n\leq 4$,  for $d=1,2,3$ in case $n=5$, for $d=1,2$ in
case $n=6$, and never happens for $n>7$
unless $d=1,$ as we quoted in the introduction.
\end{remark}

\begin{proof}
The determinant of the resultant complex provides a determinantal 
formula for  $\Rd(\efen)$ precisely when $C^{-2} = C_1 =0.$
This is respectively equivalent to the inequalities
$$ t < d_1 + d_2$$
and
$$ t_n -t = d_1 + \dots d_n - n - t  < d_1 + d_2,$$
from which the lemma follows easily.  We have decreased the
right hand side of (\ref{squareq}) by a unit in order to allow for
a natural number $t$ satisfying (\ref{squaret}).
\end{proof}
\begin{corollary}
\label{7}
For all $n \geq 7$ there exists a determinantal formula only if
$d_1=d_2=d_3=1$ and $ n-3 \leq d_4 + \dots + d_n < n,$
which forces all $d_i$ to be $1$ or at most, all of them equal $1$
except for two of them which equal $2$, or all of them equal $1$
except for one of them which equals $3$.
\end{corollary}

The proof of the corollary follows easily from
the inequality (\ref{squareq}). In any case, if a 
determinantal formula exists,  we have a determinantal
formula for $t= [t_n/2],$ as the following proposition shows. 

\begin{proposition}
\label{small}
 If a determinantal formula given by the resultant complex exists, 
then $M_{[t_n/2]}$
is square and of the smallest possible size 
$\rho(\left[\frac{t_n}{2}\right]).$
\end{proposition}

\begin{proof}
In order to prove that  $M_{[t_n/2]}$ is square, we need to check
by Lemma \ref{square} that 
\begin{equation}
\label{entera}
d_3 + \dots + d_n -n < \left[\frac{t_n}{2}\right] < d_1 + d_2 .
\end{equation}
If there exists a determinantal formula, then the inequalily
(\ref{squareq}) holds, from which it is
 straightforward to verify that  
$$ d_3 + \dots + d_n -n < \frac{t_n}{2} < d_1 + d_2 .$$
To see that  in fact (\ref{entera}) holds,  it is enough to check
that 
$$  d_3 + \dots + d_n -n  +1/2  \not=  \frac{t_n}{2} = \frac{d_1+\dots+
d_n-n}{2} .$$
But if the equality holds, we would have that
$d_3 + \dots + d_n = d_1 + d_2 + n-1,$ which is a contradiction.
According to Corollary \ref{min}, we also know that $M_{[t_n/2]}$
has the smallest possible size. 
\end{proof}

\medskip
\section{Dixon formulas}

We prove in this section that ``affine'' Dixon formulas can in fact
be recovered in this setting.
We  first recall classical Dixon formulas to compute the resultant of 
three bivariate affine polynomials of degree $d.$  We will make a slight
change of notation in what follows. The input affine polynomials (having
monomials of degree at most $d$ in two variables $(X_1,X_2)$) will be denoted
$f_1, f_2,f_3 $ and we will use capital letters $F_1,F_2,F_3$ to denote
the homogeneous polynomials in three variables given by their respective
homogenizations (with homogeneizing variable $X_3$). 
Dixon (cf. \cite{Dix}) proposed the following determinantal
formula to compute the resultant ${\rm Res}_{d,d ,d}(f_1,f_2,f_3) $
$={\rm Res}_{d,d ,d}(F_1,F_2,F_3)$:

\smallskip

Let $\B(X_1,X_2,Y_1,Y_2)$ denote the polynomial  obtained by
dividing the following determinant by $(X_1-Y_1) (X_2-Y_2)$:

$$\det \left(
\begin{array}{ccc}
f_1(X_1,X_2) &f_2(X_1,X_2) &f_3(X_1,X_2)\\
f_1(Y_1,X_2) &f_2(Y_1,X_2) &f_3(Y_1,X_2)\\
f_1(Y_1,Y_2) &f_2(Y_1,Y_2) &f_3(Y_1,Y_2).
\end{array}
\right)
$$
Note that by performing row operations we have that
$\B(X_1,X_2,Y_1,Y_2)$ equals the determinant of
the matrix
$$\det \left(
\begin{array}{ccc}
\Delta_{11}& \Delta_{21}&\Delta_{31} \\
\Delta_{12}& \Delta_{22}&\Delta_{32}\\
f_1(Y_1,Y_2) &f_2(Y_1,Y_2) &f_3(Y_1,Y_2),
\end{array}
\right)
$$
where $\Delta_{ij}$ are as in (\ref{deltaij}).
Write 
$$\B(X_1,X_2,Y_1,Y_2) = \sum_{|\beta| \leq 2d-2}
B_\beta(X_1,X_2) Y_1^{\beta_1} Y_2^{\beta_2}.$$

Set $A :=\ZZ[a],$ where $a$ denotes one indeterminate 
for each  coefficient of $f_1,f_2,f_3.$
Let $S$ denote the free module over $A$ with 
basis $\ba$ given by all monomials in two variables of
degree less or equal than $d-2,$
which has an obvious isomorphism with the free
module $S'$  over $A$ with basis ${\ba}' $ given by all monomials in
three variables of degree equal to $d-2.$ 
The monomial basis of all polynomials in two variables
of degree less or equal than $2d-2$ will be denoted by
$\ca.$
\par
Let $M$ be the square matrix of size $2 d^2-d$
whose columns are indexed by $\ca$ 
and whose rows contain consecutively the
expansion in the basis $\ca$ of $m \cdot f_1,$ 
of $m \cdot f_2,$ and of $m \cdot f_3,$ where $m$
runs in the three cases over $\ba$, and finally, the expansion
in the basis $\ca$ of all  $B_\beta, \, |\beta| \leq d-1.$
Then, Dixon's formula says that
$$ {\rm Res}_{d,d ,d}(f_1,f_2,f_3)  = \pm \det M.$$

Here, $d_1=d_2=d_3=d$ and $n=3,$ so that (\ref{squareq}) holds
and by (\ref{squaret}) there is a determinantal formula for each
$t$ such that $ d-3 \, < \, t \, < 2d.$ So, one possible choice is $t= 2d-2.$
Then, $ t_3-t= d-1 <d, $ which implies $<F_1,F_2,F_3>_{t_3-t} =0.$ 
Also, $t-d =d-2 <d,$
and therefore $S^{t,i} = S',$ for all $i=1,2,3.$ 
\par
Let $\Delta(X_1,
X_2,X_3,
Y_1,Y_2,Y_3) = \sum_{|\gamma| \leq 3d-3} \Delta_\gamma(X) Y^\gamma$
be the Bezoutian associated with the homogeneous polynomials $F_1,
F_2,F_3.$  We know that ${\rm Res}_{d,d ,d}(F_1,F_2,F_3)  = \pm \det M_{2d-2}.$
In this case, the transposed matrix $M_{2d-2}^t$ is a square matrix of the same size as $M$,
and it is obvious that their $3 d (d-1)/2$ first rows coincide (if the columns are
ordered conveniently). According to
(\ref{eq:psi1}), the last $(d+1)d/2$ rows of $M_{2d-2}^t$ contain the
expansion in the  basis ${\ba}' $ of all
$\Delta_\gamma \, , \, |\gamma| = d-1.$

\begin{proposition} 
\label{dixon} The  ``affine''  matrix $M$ and the ``homogeneous''
matrix $M_{2d-2}^t$ coincide.
\end{proposition}

\begin{proof}
Denote $P(X_1,X_2,X_3,Y_1,Y_2,t)$ the homogeneous polynomial of
degree $3d-2$ in $6$ variables  obtained by dividing the following
determinant by $(X_1-Y_1)(X_2-Y_2):$
$$
\det \left(
\begin{array}{ccc}
\Delta_{1,1}(F)& \Delta_{2,1}(F)&\Delta_{3,1}(F) \\
\Delta_{1,2}(F)& \Delta_{2,2}(F)&\Delta_{3,2}(F)\\
F_1(Y_1,Y_2,t) &F_2(Y_1,Y_2,t) &F_3(Y_1,Y_2,t)
\end{array}
\right),
$$
where
$$\Delta_{i,1}(F):= F_i(X_1,X_2,X_3) - F_i(Y_1,Y_2,X_3) , \, i=1,2,3$$
and
$$\Delta_{i,2}(F) := F_i(Y_1,X_2,X_3) -F_i(Y_1,Y_2,X_3) , \, i=1,2,3.$$

It is easy to check that 
\begin{equation}
\label{uno}
(X_3 -Y_3) \Delta(x,y) = P(X_1,X_2,X_3,Y_1,Y_2, X_3) - 
P( X_1,X_2,X_3,Y_1,Y_2, Y_3)
\end{equation}
and that
\begin{equation}
\label{dos}
P(X_1,X_2,1,Y_1,Y_2,1) = \B(X_1,X_2,Y_1,Y_2)
\end{equation}
We are looking for the elements in $\B(X_1,X_2,Y_1,Y_2)$ of degree  
less or equal than $d-1$ in the variables $Y_1,\ Y_2.$ 
But it is easy to check that
$\deg_y\left(P( X_1,X_2,X_3,Y_1,Y_2, Y_3)\right)\geq d.$ This, combined
with the
equality given in (\ref{uno}), implies that, for each $1\leq j\leq d-1:$
$$X_3\,\sum_{|\gamma|=j}\Delta_{\gamma}(X)Y^\gamma-Y_3\,
\sum_{|\gamma|=j-1}\Delta_{\gamma}(X)Y^\gamma $$
is equal to the piece of degree $j$ in the variables $Y_i$ of 
the polynomial $P(X_1,X_2,X_3,Y_1,Y_2, Y_3).$
\par
Besides, this polynomial does not depend on $Y_3,$ so the following
formula holds for every pair $\gamma,\ \tilde{\gamma}$ such that
$\gamma=\tilde{\gamma}+(0,0,k), |\gamma|=j:$
\begin{equation}
\label{clave}
X_3^k\Delta_\gamma(X) = \Delta_{\tilde{\gamma}}(X).
\end{equation}
This allows us to compute $\Delta_{\gamma}(X)$ for every $|\gamma|=d-1,$
in terms of the homogeneization of $B_{(\gamma_1,\gamma_2)}.$ 
{}From equation (\ref{clave}), the claim follows  straightforwardly.
\end{proof}

We conclude that  Dixon's formula can be
viewed as a particular case of the determinantal expressions that
we addressed. Moreover,  Proposition \ref{dixon} can be extended to any
number of variables and all Dixon matrices as in \cite[\S 3.5]{EM}
can be recovered in degrees $t$ such that
$\psi^*_{2,t_n-t} =0,$ i.e. such that $t_n \geq t > t_n - \min \{d_1,\dots, d_n\}.$
As we have seen, all one can hope in general is the explicit
quotient formula we give in Theorem \ref{mainth}.  In fact, we have
the following consequence of Lemma \ref{square}

\begin{lemma}
There exists a determinantal Dixon formula if and only if $n=2,$ or
$n=3$ and $d_1=d_2=d_3,$ i.e. in the case considered by Dixon.
\end{lemma}

\begin{proof} 
Assume $d_1 \leq d_2 \leq \dots \leq d_n.$  If
inequality (\ref{squaret}) is verified for  $t > t_n -d_1,$ we deduce
that
\begin{equation}
\label{squared}
 (n-2) d_1 -n \leq
d_3 + \dots + d_n -n < d_1 -2,
\end{equation}
and so $(n-3) d_1 < n-2.$ This equality cannot hold for any
natural number $d_1$ unless $n \leq 3.$ It is easy to check that for $n=2$
there exist a determinantal Dixon formula  for any value of $d_1, d_2.$
In case $n=3$, (\ref{squared}) implies that $d_3 < d_1 +1.$  Then, $d_1=d_2=d_3,$
as claimed.
\end{proof}

\medskip

\section{ Other known formulas and some extensions}

We can recognize other well known determinantal formulas for resultants
in this setting.

\subsection{Polynomials in one variable}

Let 
$$f_1(x) = \sum_{j=0}^{d_1} a_j x^j \, \,  , \, \, 
f_2(x) = \sum_{j=0}^{d_2} b_j x^j $$ be generic univariate polynomials
(or their homogenizations in two  variables) of degrees $d_1 \leq d_2.$
In this case,  inequality (\ref{squaret}) is verified
for all $ t  = 0, \dots, d_1 +d_2 -1$ and so we have a determinantal
formula for all such $t$. Here, $t_2 = d_1 + d_2 -2.$ When
$t= d_1+d_2 -1 = t_2 +1$ we have the classical Sylvester
formula.

Assume $d_1 = d_2= d$  and write
$$ \frac { f_1(x) f_2(y) - f_1(y) f_2(x)} {x-y} = \sum_{i,j=0}^d
c_{ij} x^i y^j.$$
Then, the classical  {\it B\'ezout \/} formula for the
 resultant between $f_1$ and $f_2$ says that
$${\rm Res}_{d,d}(f_1,f_2) = \det  (c_{ij}).$$
It is easy to see that we obtain precisely this formulation for
$t= d-1.$ For other values of $t$ we get formulas
interpolating between Sylvester and B\'ezout as in
\cite[Ch. 12]{gkz}, even in case $d_1 \neq d_2$.  It is
easy to check that the smallest possible matrix has size $d_2.$ 

Suppose for example that $d_1 =1, d_2 = 2.$ In this case,
$\left[ t_2/2\right] = \left[ 1/2 \right] = 0,$ and $M_0$ is
a $2\times 2 $ matrix representing a map from
$S_1^*$ to $S_0 \oplus S_0^*,$ whose determinant
equals the resultant
$${\rm Res}_{1,2}(f_1,f_2) = a_1^2 b_0 - a_0 a_1 b_1 + b_2 a_0^2.$$
If we write $f_1(x) = 0 x^2 + a_1 x + a_0$ and we use the classical
Bezout formula for $d=2,$ we would also get a $2 \times 2$
matrix but whose determinant equals $b_2 \cdot {\rm Res}_{1,2}(f_1,f_2).$
The exponent $1$ in $b_2$ is precisely the difference $d_2 - d_1.$

\subsection{ Sylvester formula for three ternary quadrics}

Suppose that $n=3, \ d_1=d_2=d_3=2$  and $2 \not=0.$ Let $J$ denote
the Jacobian determinant associated with the homogeneous
polynomials $f_1,f_2,f_3.$  A beautiful classical formula due
to Sylvester says that the resultant ${\rm Res}_{2,2 ,2}(f_1,f_2,f_3)$
can be obtained as $1/512$ times the determinant of the
$6 \times 6 $ matrix whose columns are indexed by the
monomials in $3$ variables of degree $2$ and whose rows
correspond to the expansion in this monomial basis of
$f_1, f_2, f_3, \frac{\partial J}{\partial X_1}, \frac{\partial J}{\partial X_2}$ and
$\frac{\partial J}{\partial X_3}.$ In this case, $\left[ t_3/2\right]=
\left[ 3/2\right] = 1,$  and by Lemma \ref{square} we have a determinantal formula in this degree
 since $ 2 - 3 < 1 < 4.$ From Euler equations
$$ 2 f_i(X) = \sum_{j=1}^3 X_j  \frac{\partial f_i(X)}{\partial X_j},$$
we can write
$$ \begin{array}{ccl}
2\left(f_i(X)-f_i(Y)\right)&=&\sum_{j=1}^3 {\left(
X_j  \frac{\partial f_i(X)}{\partial X_j} - Y_j  \frac{\partial f_i(Y)}
{\partial Y_j}\right)}\\
&=& \sum_{j=1}^3 {(X_j-Y_j)\frac{\partial f_i(X)}{\partial X_j}+
Y_j\left(\frac{\partial f_i(X)}{\partial X_j}-\frac{\partial f_i(Y)}
{\partial Y_j}\right)}\\
&=& \sum_{j=1}^3 \left((X_j-Y_j)\,\frac{\partial f_i(X)}{\partial X_j}+ 
Y_j\sum_{l=1}^3 {\frac{\partial^2 f_i(X)}{\partial X_j\partial X_l}
(X_l-Y_l)}\right).\\
\end{array}$$
Because of (\ref{deltamod}), we can compute the Bezoutian using  
$$\Delta_{ij}(X,Y):=\frac12\left( \frac{\partial f_i(X)}{\partial X_j} +
\sum_{l=1}^3{\frac{\partial^2 f_i(X)}{\partial X_l\partial X_j}} Y_l\right).$$
Using this formulation, it is not difficult to see that we can recover 
Sylvester formula from the equality $ {\rm Res}_{2,2 ,2}(f_1,f_2,f_3) = 
\det M_1.$

\subsection{Jacobian formulations}

When $t=t_n,$ one has $H_d(t)=1,$  and via  the canonical identification of  $S_0^*$ 
 with $A,$ the complex (\ref{def:complex}) reduces to the following modified Koszul Complex:
\begin{equation}
\label{jac1}
 0 \longrightarrow \, K(t)^{-n} \, 
\mathop{\longrightarrow}^{\delta_{-(n-1)}}
\, \dots \, \mathop{\longrightarrow}^{\delta_{-1}} \, A\oplus K(t)^{-1}\, 
\mathop{\longrightarrow}^{\delta_{0}} \, K(t)^{0}\,{\longrightarrow}0,
 \end{equation}
where  $\delta_0$ is the following map:
$$
\begin{array}{ccccc}
 A& \oplus& 
S_{t_n-d_1}\oplus \cdots \oplus S_{t_n-d_n} 
&\rightarrow & S_{t_n} 
\\
(\lambda&,&g_1,\ldots,g_n)&\mapsto&\lambda\,\Delta_0+\sum_{i=1}^n g_i\,f_i,\\
\end{array}
$$
and $\Delta_0:=\Delta(X,0).$ 
As a corollary of Theorem \ref{rescth} we get that, for every specialization of the coefficients, 
$\Delta_0$ is a non-zero element of the quotient if the resultant does 
not vanish.
\par
Assume that the characteristic of $k$ does not divide the product $d_1\dots d_n.$
It is a well-known fact that  the  jacobian determinant $J$ of the sequence $(\efen)$ is another 
non-zero element of degree $t_n,$ which is  a non-zero element of the quotient whenever
the resultant does not vanish (cf. for instance \cite{ss}).
In fact, one can easily check that
\begin{equation}
\label{id1}
J \ =\ d_1\dots d_n \, \Delta_0 \ \mod \left<f_1,\ldots,f_n\right>.
\end{equation}
In \cite{CDS}, the same complex  
is considered in a more general toric setting, but  using $J$ instead of $\Delta_0$ . 
Because of  (\ref{id1}), we can recover their results in the homogeneous case.

\begin{theorem}
Consider the modified complex (\ref{jac1}) with $J$ instead of $\Delta_0.$ Then, for every 
specialization of the coefficients,
the complex is exact if and only if the resultant does not vanish. Moreover, the determinant
of the complex equals $d_1\dots d_n \, \Rd(\efen).$
\end{theorem}
We can also replace $\Delta_0$ by $J$ in Macaulay's Formula (Theorem \ref{mainth}),
and  have the following result:
\begin{theorem}
Consider the square submatrix $\tilde M_{t_n}$ which is extracted from the matrix of
$\delta_0$ in the monomial bases, choosing the same rows and columns of $M_{t_n}.$
Then, $\det(\tilde M_{t_n})\neq0,$ and we have the following formula 
\`a la Macaulay:
$$d_1\dots d_n\,\Rd(\efen) = \frac{\det(\tilde M_{t_n})}{\det(E_{t_n})}.$$
\end{theorem}

\medskip

We end the paper by addressing two natural questions that arise: 

\subsection{Different choices of monomial bases}

Following Macaulay's original ideas, one can show that
there is some flexibility in the choice of the monomial bases 
defining $S^{t,i}$  in order to get  other non-zero minors, 
of $\tilde\Psi_t$, i.e  different square matrices
$M'_t$  whose determinants are  non-zero multiples
of $\Rd(\efen)$ with different extraneous factors 
$\det(E'_t), \det(E'_{t_n-t})$, for 
appropiate square submatrices $E'_t, E'_{t_n-t}$ of $M'_t.$
Besides the obvious choices coming from a permutation
in the indices of the variables, other choices can be made as follows. 
\par
For any $i=1,\dots,n,$
set $\hat{d}_i:= (d_1, \dots,
d_{i-1},d_{i+1},\dots,d_n)$ and define $H_{\hat{d}_i}(t)$ for any positive
integer $t$ by the equality
$$\frac{\prod_{j\neq i}\left( 1-Y^{d_{j}}\right) }{\left( 1-Y\right) ^{n-1}}=
\sum_{t=0}^\infty H_{\hat{d}_i}(t).Y^t.$$
For each $t\in\NN_0,$ set  also
$\Lambda_t:=\{X^\gamma \in S_t: \ \gamma_j<d_j, \, j=1,\dots,n\}.$

We then have the following result:

\begin{proposition}
Let $M'_t$ a square submatrix of $\tilde\Psi_t$ of size $\rho(t).$
Denote its blocks as in (\ref{prima}).
Suppose that, for each $i=1,\dots,n$, the block $F_t$ 
has  {\it exactly} $H_{\hat{d}_i}(t-d_i)$ of its columns corresponding to $f_i$ 
in common with  the matrix  $D_t$  defined in (\ref{prev}) 
and, also, the block $F_{t_n-t}$ shares exactly $H_{\hat{d}_i}(t_n-t-d_i)$ 
columns
corresponding  to $f_i$  with $D_{t_n-t}.$ Then, if $\det\left(
M'_t\right)$ is not identically zero,  the
resultant $\Rd(\efen)$ can be computed as the ratio $\frac{\det( M'_t)}
{\det(\EE_t')},$ where $\EE_t'$ is made by joining two submatrices
$E'_t$ of $F_t$ and $E'_{t_n-t}$ of
$F_{t_n-t}.$ These submatrices are obtained by omitting the columns in common with
$D_t$ (resp. $D_{t_n-t}$) and the rows indexed by all common monomials in 
 $D_t$ (resp. $D_{t_n-t}$) and all monomials in 
$\Lambda_t$ (resp. $\Lambda_{t_n-t}$). 
\end{proposition}
We omit the proof which is rather technical, and based in
 \cite[6a]{Mac}, and \cite{Ch1},\cite{cha}.

\subsection{Zeroes at infinity}

Given a non-homogeneous system of  polynomial equations $\tilde{f}_1,
\dots,\tilde{f}_n$ in $n-1$ variables with respective degrees $d_1,\dots,d_n,$ we can
homogenize these polynomials and consider the resultant $\Rd(\efen)$
associated with their respective homogenizations $\efen.$ However,
this resultant may vanish due to common zeros of $\efen$ at infinity
in projective space
$\CP^{n-1}$ even when there is no affine common root to 
$\tilde{f}_1,= \dots =\tilde{f}_n =0.$
We can in this case extend Canny's construction \cite{C2}
of the  Generalised Characteristic Polynomial (GCP) for classical
Macaulay's matrices to the matrices $M_t$ for any natural number
$t.$  In fact, 
when we specialize  $f_i$ to $ X_i^{d_i}$ for
all $i=1,\dots,n,$ the Bezoutian is given by 
$$ \sum_{j_1=0}^{d_1-1} \cdots \sum_{j_n=0}^{d_n-1} X_1^{d_1-1-j_1}
\cdots X_n^{d_n-1-j_n} \ Y_1^{j_1} \cdots Y_n^{j_n},$$
and the specialized matrix $M_{t}(e)$ of $M_t$ has a 
single non zero entry on each row and column 
which is equal to $1,$ so that $\det(M_t(e)) = \pm 1.$  We order the
columns in such a way that $M_t(e)$ is the identity matrix.
With this convention,  define the polynomial $C_t(s)$   
by
$$C_t(s) := \frac { {\mbox {Charpoly }} (M_t) (s)} { {\mbox {Charpoly }} (\EE_t) (s)},$$
where $s$ denotes a new variable and Charpoly  means characteristic
polynomial.
We then have by the previous observation that
$$C_t(s) = \Rd(f_1 - s \ x_1^{d_1}, \dots, f_n - s \ x_1^{d_n}).$$

Moreover, this implies that $C_t(s)$ coincides with Canny's GCP
$C(s),$  but involves matrices of smaller size. 
Canny's considerations on how to compute more efficiently
the GCP also hold in this case.
Of course, it is in general much better to  find a way to
construct ``tailored'' residual 
resultants for polynomials with a generic structure which is not
dense, as in the case of sparse polynomial systems (\cite{EM}, 
\cite{gkz}).

\medskip

\noindent {\bf Acknowledgements:} 
We are grateful to J. Fern\'andez Bonder, G. Massaccesi and J.M. Rojas for
helpful suggestions. We are also grateful to David Cox for his
thorough reading of the manuscript.

\noindent{\bf Author's addresses:}

\noindent{Departamento  de Matem\'atica, F.C.E y N., UBA,
(1428) Buenos Aires, Argentina. }

\noindent{{\tt cdandrea@dm.uba.ar \hskip4.5truecm alidick@dm.uba.ar}}

\end{document}